\documentstyle[11pt,a4wide]{article}
\font \tenmsb=msbm10 scaled \magstep 1
\font \sevenmsb=msbm7 scaled \magstep 1
\font \fivemsb=msbm5 scaled \magstep 1
\newfam \msbfam
\textfont \msbfam=\tenmsb
\scriptfont \msbfam=\sevenmsb
\scriptscriptfont \msbfam=\fivemsb
\def \Bbb#1{\fam \msbfam \relax#1}

\font \teneufm=eufm10 scaled \magstep 1
\font \seveneufm=eufm7 scaled \magstep 1
\font \fiveeufm=eufm5 scaled \magstep 1
\newfam \eufmfam
\textfont \eufmfam=\teneufm
\scriptfont \eufmfam=\seveneufm
\scriptscriptfont \eufmfam=\fiveeufm
\def \frak#1{{\fam \eufmfam \relax#1}}

\title{\bf ON FUNCTION THEORY IN QUANTUM DISC: INVARIANT KERNELS}

\author{\sl D. Shklyarov \and \sl S. Sinel'shchikov
\and \sl L. Vaksman \thanks{Partially supported by the grant
INTAS-94-4720}}

\date{\tt Institute for Low Temperature Physics \& Engineering\\
National Academy of Sciences of Ukraine}

\newtheorem{theorem}{Theorem}[section]
\newtheorem{lemma}[theorem]{Lemma}
\newtheorem{proposition}[theorem]{Proposition}
\newtheorem{corollary}[theorem]{Corollary}

\begin{document}

\maketitle

 In our earlier work \cite{SSV1} some results on integral representations of
functions in the quantum disc were formulated without proofs. It was then
shown in \cite{SSV2} that the validity of those results is related to the
invariance of kernels of integral operators considered in \cite{SSV1}. We
introduce in this work a method which allows us to prove the invariance of
the above kernels.

 Note that the invariant kernels under consideration may be treated as
generating functions from the viewpoint of the theory of basic
hypergeometric series \cite{GR}. A simplest example of relationship between
generating functions and q-special functions presents the q-binomial of
Newton: let $ab=qba$, then
$$(a+b)^n=\sum_{k=0}^n \frac{(q;q)_n}{(q;q)_k(q;q)_{n-k}}b^ka^{n-k},\qquad
(q;q)_m=(1-q)(1-q^2)\ldots(1-q^m).$$
Another example is the q-analogue of Van der Vaerden generating function for
Clebsch-Gordan coefficients \cite{V}.

\bigskip \setcounter{equation}{0}

\section{The principal homogeneous space $\widetilde{X}$}

 Consider the real Lie group
$$SU(1,1)=\left \{\left.\pmatrix{t_{11} & t_{12}\cr t_{21} & t_{22}}\in
SL_2({\Bbb
C})\right|\overline{t}_{11}=t_{22},\:\overline{t}_{12}=t_{21}\right \}$$
and the element $w=\pmatrix{0 & -1 \cr 1 & 0}\in SL_2({\Bbb C})$. Evidently,
$w \notin SU(1,1)$.

 The subgroup $SU(1,1)$ acts on $SL_2({\Bbb C})$ via right shifts, and
$\widetilde{X}=w^{-1}\cdot SU(1,1)$ is an orbit of this action. Obviously,
$$\widetilde{X}=\left \{\left.\pmatrix{t_{11} & t_{12}\cr t_{21} &
t_{22}}\in SL_2({\Bbb
C})\right|\overline{t}_{11}=-t_{22},\:\overline{t}_{12}=-t_{21}\right \}.$$

 Now turn to the construction of a q-analogue for the homogeneous space
$\widetilde{X}$.

 Remind first the definition of the algebra ${\Bbb C}[SL_2]_q$ of regular
functions on the quantum group $SL_2$ (see \cite{CP}). This algebra is
determined by its generators $\{t_{ij}\}$, $i,j=1,2$, and the relations:
\begin{equation}\label{tijcr}
\left \{\begin{array}{lcl}t_{11}t_{12}=qt_{12}t_{11}&,&
t_{21}t_{22}=qt_{22}t_{21},\\ t_{11}t_{21}=qt_{21}t_{11}&,&
t_{12}t_{22}=qt_{22}t_{12},\\ t_{12}t_{21}=t_{21}t_{12}&,&
t_{11}t_{22}-t_{22}t_{11}=(q-q^{-1})t_{12}t_{21}\\
t_{11}t_{22}-qt_{12}t_{21}=1&.&\end{array}\right.
\end{equation}

 Equip ${\Bbb C}[SL_2]_q$ with the involution given by
\begin{equation}\label{inv}
t_{11}^*=-t_{22},\qquad t_{12}^*=-qt_{21}.
\end{equation}
We thus get an involutive algebra which is denoted by ${\rm
Pol}(\widetilde{X})_q$.

 The simplest argument in favor of our choice of involution is that at the
limit $q \to 1$ one gets the system of equations $t_{11}^*=-t_{22}$,
$t_{12}^*=-t_{21}$, which distinguishes the orbit $\widetilde{X}$.

 A more complete list of arguments and the "algorithm" of producing
involutions which leads to the $*$-algebra ${\rm Pol}(\widetilde{X})_q$ are
expounded in the appendix.

 It is well known \cite{CP} that $w^{-1}=\pmatrix{0 & 1 \cr -1 & 0}$ admits
a quantum analogue $w_q^{-1}\in{\Bbb C}[SL_2]_q^*$. One can deduce from the
results of appendix, in particular, that this linear functional is real.
Besides that, it is shown in appendix that ${\rm Pol}(\widetilde{X})_q$ is a
covariant $*$-algebra.

 Note that the structure of a $U_q \frak{su}(1,1)$-module algebra in ${\rm
Pol}(\widetilde{X})_q$ is determined by the relation (1.1) of \cite{SSV2}
and the relations
\begin{equation}\label{pol-mod}
\left \{\begin{array}{rcl}\pmatrix{X^+t_{11}& X^+t_{12}\cr X^+t_{21}&
X^+t_{22}}&=&\pmatrix{t_{11}& t_{12}\cr t_{21}& t_{22}}\pmatrix{0 & 1 \cr 0
& 0}\\ \pmatrix{X^-t_{11}& X^-t_{12}\cr X^-t_{21}&
X^-t_{22}}&=&\pmatrix{t_{11}& t_{12}\cr t_{21}& t_{22}}\pmatrix{0 & 0 \cr 1
& 0}\\ \pmatrix{Ht_{11}& Ht_{12}\cr Ht_{21}& Ht_{22}}&=&\pmatrix{t_{11}&
t_{12}\cr t_{21}& t_{22}}\pmatrix{1 & 0 \cr 0 & -1}\end{array}\right..
\end{equation}

\bigskip \setcounter{equation}{0}

\section{A completion}

 Turn to a construction of the algebra $D(\widetilde{X})_q$ of finite
functions and the bimodule $D(\widetilde{X})_q^\prime$ of distributions on
the quantum principal homogeneous space $\widetilde{X}$.

 Let
$$M^\prime=\{(i_1,i_2,j_1,j_2)\in{\Bbb Z}_+^4|\;i_2 \cdot j_2=0 \},\qquad
M^{\prime \prime}=\{(i_1,i_2,j_1,j_2)\in{\Bbb Z}_+^4|\;i_1j_1=i_2j_2=0
\}.$$

 It follows from the definition of ${\rm Pol}(\widetilde{X})_q$ that every
element $f \in{\rm Pol}(\widetilde{X})_q$ admits a unique
decomposition \begin{equation}\label{exp1}
f=\sum_{(i_1,i_2,j_1,j_2)\in M^\prime}a_{i_1i_2j_1j_2}(f)
t_{11}^{i_1}t_{12}^{i_2}t_{12}^{*j_2}t_{11}^{*j_1}.
\end{equation}

 Let $x=t_{12}t_{12}^*$. By a virtue of (\ref{exp1}), each element $f
\in{\rm Pol}(\widetilde{X})_q$ admits a unique decomposition
\begin{equation}\label{exp2}
f=\sum_{(i_1,i_2,j_1,j_2)\in M^{\prime
\prime}}t_{11}^{i_1}t_{12}^{i_2}\psi_{i_1i_2j_1j_2}(x)
t_{12}^{*j_2}t_{11}^{*j_1}.
\end{equation}

\medskip

 {\sc Remark 2.1.} In view of the definitions to be imposed below, it
should be noted that the sums in (\ref{exp1}), (\ref{exp2}) are finite, the
coefficients $\psi_{i_1i_2j_1j_2}(x)$ are polynomials, and the functionals
$l_{i_1i_2j_1j_2}^\prime:f \mapsto a_{i_1i_2j_1j_2}(f)$,
$(i_1,i_2,j_1,j_2)\in M^\prime$, $l_{i_1i_2j_1j_2}^{\prime \prime}:f \mapsto
\psi_{i_1i_2j_1j_2}(q^{-2m})$, $(i_1,i_2,j_1,j_2)\in M^{\prime \prime}$, $m
\in{\Bbb Z}_+$, are linear.

\medskip \stepcounter{theorem}

 Equip the vector space ${\rm Pol}(\widetilde{X})_q$ with the weakest
topology in which all the functionals $l_{i_1i_2j_1j_2}^{\prime \prime}$,
$(i_1,i_2,j_1,j_2)\in M^{\prime \prime}$, $m \in{\Bbb Z}_+$, are continuous.

 A completion of the above Hausdorff topological vector space will be called
the space of distributions and denoted by $D(\widetilde{X})_q^\prime$.

\medskip

 {\sc Remark 2.2.} The vector space ${\rm Pol}(\widetilde{X})_q$ is equipped
with the structure of algebra, and hence with the structure of ${\rm
Pol}(\widetilde{X})_q$-bimodule. This structure is extendable by a
continuity onto the completion $D(\widetilde{X})_q^\prime \supset{\rm
Pol}(\widetilde{X})_q$.

\medskip \stepcounter{theorem}

 $D(\widetilde{X})_q^\prime$ will be identified with the space of formal
series (\ref{exp2}), whose coefficients are functions on $q^{-2{\Bbb Z}_+}$.
The topology in the space of such series is that of pointwise convergence of
the coefficients $\psi_{i_1i_2j_1j_2}(x)$.

 A distribution $f \in D(\widetilde{X})_q^\prime$ is said to be finite if
$\#\{(i_1,i_2,j_1,j_2,m)|\;\psi_{i_1i_2j_1j_2}(q^{-2m})\ne 0 \}<\infty$. The
vector space of finite functions will be denoted by $D(\widetilde{X})_q$.
(It is easy to present a non-degenerate pairing $D(\widetilde{X})_q^\prime
\times D(\widetilde{X})_q \to{\Bbb C}$ which establishes an isomorphism
between $D(\widetilde{X})_q^\prime$ and a vector space dual to
$D(\widetilde{X})_q$).

\medskip

\begin{proposition} The structure of a covariant {\tt $*$-algebra} is
transferred by a continuity from ${\rm Pol}(\widetilde{X})_q$ onto
$D(\widetilde{X})_q$. The structure of a covariant
{\tt $D(\widetilde{X})_q$-bimodule} is transferred by a continuity from
$D(\widetilde{X})_q$ onto $D(\widetilde{X})_q^\prime$.
\end{proposition}

\smallskip

 {\bf Proof.} A verification of the covariance of the algebra
$D(\widetilde{X})_q$ and the bimodule $D(\widetilde{X})_q^\prime$ reduces to
the proof of the identities whose validity is already known in the case when
$\psi_{i_1i_2j_1j_2}(x)$ are polynomials. Now what remains is to observe
that the polynomials are dense in the space of functions $\psi(x)$
equipped with the topology of pointwise convergence. \hfill $\Box$

\medskip

\begin{lemma} The formal series (\ref{exp1}) presents a finite function on
the quantum principal homogeneous space iff the following two conditions are
satisfied:
\\ i) $\#\{(i_2,j_2)\in{\Bbb Z}_+^2|\;i_2j_2=0,\,\exists(i_1,j_1)
\in{\Bbb Z}_+^2:\:a_{i_1i_2j_1j_2}(f)\ne 0 \}<\infty$,
\\ ii) $\exists N \in{\Bbb N}:\quad f \cdot t_{11}^N=t_{11}^{*N}\cdot f=0$.
\end{lemma}

\medskip

\begin{proposition} There exists a unique nonzero element $e_0 \in
D(\widetilde{X})_q$ such that
\begin{equation}\label{e0i}
t_{12}e_0=e_0t_{12},\qquad t_{12}^*e_0=e_0t_{12}^*,
\end{equation}
\begin{equation}\label{e0ii}t_{11}^*e_0=e_0t_{11}=0,\end{equation}
\begin{equation}\label{e0iii}e_0 \cdot e_0=e_0.\end{equation}
\end{proposition}

\smallskip

 {\bf Proof.} $e_0$ satisfies (\ref{e0i}) iff $e_0=\displaystyle
\sum_{i_2j_2=0}t_{12}^{i_2}\psi_{0i_20j_2}(x)t_{12}^{*j_2}$. (\ref{e0iii})
means that
\begin{equation}\label{psi0i}e_0=\psi_{0000}(x),\end{equation}
with $(\psi_{0000}(x))^2=\psi_{0000}(x)$. Finally, (\ref{e0ii}) is satisfied
iff
\begin{equation}\label{psi0ii}
\psi_{0000}(x)=\left \{\begin{array}{ccl}1 &,& x=1 \\ 0 &,& x \in
\{q^{-2},q^{-4},q^{-6},\ldots \}\end{array}\right..
\end{equation}
Thus $e_0$ exists, is unique and determined by (\ref{psi0i}),
(\ref{psi0ii}). \hfill $\Box$

\medskip

 It is easy to prove that the structure of a covariant $*$-algebra can be
transferred by a continuity from ${\rm Pol}(\widetilde{X})_q$ onto ${\rm
Fun}(\widetilde{X})_q \stackrel{\rm def}{=}{\rm
Pol}(\widetilde{X})_q+D(\widetilde{X})_q$. The structure of a $*$-algebra
is given by (\ref{e0i}) -- (\ref{e0iii}) and $e_0^*=e_0$, and the $U_q
\frak{su}(1,1)$-action by
\begin{equation}\label{e0iv}
He_0=0,\qquad X^+e_0=c_+t_{11}e_0t_{12}^*,\qquad
X^-e_0=c_-t_{12}e_0t_{11}^*,
\end{equation}
with $c_-=-{\textstyle e^{-3h/4}\over \textstyle 1-e^{-h}}$,
$c_+=-{\textstyle e^{-5h/4}\over \textstyle 1-e^{-h}}$. The proof of
(\ref{e0iv}) is just the same as that of similar relations (3.5) in
\cite{SSV2}.

 Finally note that for any $f \in D(U)_q$ there exists a unique
decomposition \begin{equation} f=\sum_{(i_1,i_2,j_1,j_2)\in
M^\prime}c_{i_1i_2j_1j_2}
t_{11}^{i_1}t_{12}^{i_2}e_0t_{12}^{*j_2}t_{11}^{*j_1},
\end{equation}
with $c_{i_1i_2j_1j_2}\in{\Bbb C}$ and only finitely many of those being
non-zero.

\bigskip \setcounter{equation}{0}

\section{Quantum homogeneous space $X$}

 In the case $q=1$ the orbit $\widetilde{X}=w^{-1}\cdot SU(1,1)$ can be
equipped with a structure of a homogeneous $SU(1,1)\times SU(1,1)$-space:
$$(g_1,g_2):x \mapsto (w^{-1}g_1w)xg_2^{-1}.$$
In particular, the action of the one-parameter subgroup
$$i_ \varphi:\pmatrix{t_{11} & t_{12}\cr t_{21} & t_{22}}\mapsto
g_1(\varphi)\cdot \pmatrix{t_{11} & t_{12}\cr t_{21} & t_{22}},\qquad
g_1(\varphi)=\pmatrix{e^{i \varphi} & 0 \cr 0 & e^{-i \varphi}}$$
commutes with the right multiplication by $g_2^{-1}\in SU(1,1)$.

 These constructions can be transferred onto the case $0<q<1$ via (A.1).
Here we restrict ourselves to introducing the one-parameter group $i_
\varphi$, $\varphi \in{\Bbb R}/(2 \pi{\Bbb Z})$, of automorphisms of the
covariant $*$-algebra ${\rm Fun}(\widetilde{X})_q={\rm
Pol}(\widetilde{X})_q+D(\widetilde{X})_q$:
$$\pmatrix{i_ \varphi(t_{11}) & i_ \varphi(t_{12}) \cr i_ \varphi(t_{21}) &
i_ \varphi(t_{22})}=\pmatrix{e^{i \varphi} & 0 \cr 0 & e^{-i
\varphi}}\pmatrix{t_{11} & t_{12}\cr t_{21} & t_{22}},\qquad i_
\varphi(e_0)=e_0.$$
Evidently, $i_ \varphi \cdot{\rm Pol}(\widetilde{X})_q={\rm
Pol}(\widetilde{X})_q$, $i_ \varphi \cdot
D(\widetilde{X})_q=D(\widetilde{X})_q$, and the operators $i_ \varphi$ are
extendable by a continuity onto the entire space
$D(\widetilde{X})_q^\prime$.

 In the case $q=1$ the orbits of the one-parameter transformation group $i_
\varphi:\widetilde{X}\to \widetilde{X}$ form a homogeneous space $X$ of the
group $SU(1,1)$ isomorphic to the unit disc $U \subset{\Bbb C}$. We intend to
produce a similar isomorphism in the quantum case $0<q<1$. Let us start
with introducing the spaces of invariants of the one-parameter group $i_
\varphi$:
$$D(X)_q=\{f \in D(\widetilde{X})_q|\;i_ \varphi(f)=f,\:\varphi \in{\Bbb
R}/(2 \pi{\Bbb Z})\},$$
$$D(X)_q^\prime=\{f \in D(\widetilde{X})_q^\prime|\;i_
\varphi(f)=f,\:\varphi \in{\Bbb R}/(2 \pi{\Bbb Z})\}.$$
The elements of the covariant $*$-algebra $D(X)_q$ (the covariant
$D(X)_q$-bimodule $D(X)_q^\prime$) will be called finite functions (resp.
distributions) on the quantum homogeneous space $X$.

 In our work \cite{SSV2} a finite function in the quantum disc $f_0$ such
that $z^*f_0=f_0z=0$, $f_0 \cdot f_0=f_0$ is introduced.

\medskip

\begin{proposition} There exists a unique isomorphism $i:D(U)_q {{\atop
\displaystyle \to} \atop{\displaystyle \sim \atop}}D(X)_q$ of covariant
$*$-algebras such that $i:f_0 \mapsto e_0$.  The operator $i$ admits an
extension up to an isomorphism $\overline{i}:D(U)_q^\prime{{\atop
\displaystyle \to}\atop{\displaystyle \sim \atop}}D(X)_q^\prime$ of $U_q
\frak{su}(1,1)$-modules by a continuity, which agrees with the bimodule
structures.
\end{proposition}

\smallskip

 {\bf Proof.} The uniqueness of $i$ is obvious since $f_0$
generates $U_q \frak{su}(1,1)$-module $D(U)_q$ (see \cite{SSV2}).

 Prove the existence of $i$.

 Consider the linear subspace $F$ of those $f \in D(\widetilde{X})_q^\prime$
for which only finitely many terms in the formal series (\ref{exp2}) are
non-zero.  It is easy to show that the structure of a covariant $*$-algebra
is extendable by a continuity from ${\rm Pol}(\widetilde{X})_q$ onto $F$. It
is worthwhile to note that $e_0 \in F$ and the elements $t_{12}$, $t_{12}^*$
are invertible in the algebra $F$:
$$t_{12}^{-1}=x^{-1}\cdot t_{12}^*,\qquad t_{12}^{*-1}=t_{12}\cdot x^{-1}.$$
In \cite{SSV2} one can find a description of the covariant $*$-algebra ${\rm
Fun}(U)_q={\rm Pol}({\Bbb C})_q+D(U)_q$ in terms of generators and
relations. It follows from that description and the results of the previous
section that the map
\begin{equation}\label{i}
i:z \mapsto qt_{11}t_{12}^{-1},\qquad i:z^*\mapsto t_{21}^{-1}t_{22},\qquad
i:f_0 \mapsto e_0
\end{equation}
is extendable up to a homomorphism of covariant $*$-algebras $i:{\rm
Pol}(U)_q \to F$. It remains to note that $D(U)_q \subset{\rm Fun}(U)_q$.
The possibility of extending by a continuity of $i$ as well as the
properties of its extension $\overline{i}$ are easily derivable from
(\ref{exp2}). \hfill $\Box$

\medskip

 We shall not distinguish in the sequel the functions in the quantum disc
and their images under the embedding $i$.

 Let $F_m=\{f \in F|\;i_ \varphi(f)=e^{im \varphi}f,\:\varphi \in{\Bbb R}/(2
\pi{\Bbb Z})\}$. Note that for all $m \in{\Bbb Z}$
\begin{equation}
F_m \cap D(\widetilde{X})_q=D(U)_q \cdot t_{12}^m.
\end{equation}

 Finally, we present an explicit form of an invariant integral on the
quantum principal homogeneous space. Consider the linear functional
$\nu:D(\widetilde{X})_q \to{\Bbb C}$,
\begin{equation}\label{nu}
\int \limits_{\widetilde{X}_q}fd \nu=(1-q^2)\sum_{m=0}^\infty
\psi_{0000}(q^{-2m})q^{-2m},
\end{equation}
determined by the coefficient $\psi_{0000}$ in the expansion (\ref{exp2}).

\medskip

\begin{proposition} The linear functional (\ref{nu}) is an invariant
integral.
\end{proposition}

\smallskip

 {\bf Proof.} It is well known (see \cite{SSV2}) that the functional $D(U)_q
\to{\Bbb C}$,
$$\sum_{j>0}z^j \psi_j(y)+\psi_0(y)+
\sum_{j>0}\psi_{-j}(y)z^{*j}\mapsto(1-q^2)\sum_{j=0}^\infty
\psi_{0000}(q^{2j})q^{-2j}$$
is an invariant integral (here $y=1-zz^*$). Hence, by a virtue of
proposition 3.1 and (\ref{i}), the restriction of the linear functional
(\ref{nu}) onto the covariant $*$-subalgebra $D(U)_q \simeq D(X)_q$ is an
invariant integral. It remains to elaborate the fact that the averaging
operator
$$j:D(\widetilde{X})_q \to D(X)_q;\qquad j:f \mapsto{1 \over 2 \pi}\int
\limits_0^{2 \pi}i_ \varphi(f)d \varphi$$
is a morphism of $U_q \frak{su}(1,1)$-modules. \hfill $\Box$

\bigskip \setcounter{equation}{0}

\section{Asymptotic cones $\widetilde{\Xi}$ and $\Xi$}

 The definition of the covariant $*$-algebra ${\rm Pol}(\widetilde{X})_q$
involves the generators $t_{ij}$, $i,j=1,2$, and the relations (\ref{tijcr}),
(\ref{inv}), (\ref{pol-mod}). All those relations are homogeneous in
$t_{ij}$ except the relation
\begin{equation}\label{n-hom}t_{11}t_{22}-qt_{12}t_{21}=1.\end{equation}
Hence the substitution
\begin{equation}\label{subst}
t_{ij}=q^{-N}\cdot t_{ij}^{(N)},\qquad N \in{\Bbb N},
\end{equation}
together with the succeeding passage to a limit as $N \to +\infty$ lead to
the same system of relations except (\ref{n-hom}).  The latter relation at
the limit $N \to +\infty$ changes to
$$t_{11}t_{22}-qt_{12}t_{21}=0.$$
The covariant $*$-algebra given by the initial list of generators and the
new list of relations (cf. (\ref{tijcr}), (\ref{inv}), (\ref{pol-mod})) is
denoted by ${\rm Pol}(\widetilde{\Xi})_q$. It is a q-analogue of the
polynomial algebra on the cone $\widetilde{\Xi}=\{(t_{11},t_{12})\in{\Bbb
C}^2|\;|t_{11}|=|t_{12}|\}$.

 The automorphisms $i_ \varphi$ are defined in the same way as in section 3.
Their invariants constitute a covariant $*$-algebra which will be denoted in
the sequel by ${\rm Pol}(\Xi)_q$. This is a q-analogue of the polynomial
algebra on the cone
$$\Xi=\{(t_{11},t_{12})|\;t_{11}\in{\Bbb C},\:t_{12}\in{\Bbb
R},\:t_{12}=|t_{11}|\}.$$

 Let us pass from polynomials to distributions. As above, set up
$x=t_{12}t_{12}^*$. Each element $f \in{\rm Pol}(\widetilde{\Xi})_q$ admits
a unique decomposition
\begin{equation}\label{exp3}
f=\sum_{(i_1,i_2,j_1,j_2)\in M^{\prime \prime}}t_{11}^{i_1}t_{12}^{i_2}
\psi_{i_1i_2j_1j_2}(x)t_{12}^{*j_2}t_{11}^{*j_1}.
\end{equation}
Remind (see section 2) that in the case of ${\rm Pol}(\widetilde{X})_q$ a
similar decomposition it possible to impose a topology. Specifically, we
equipped the vector space ${\rm Pol}(\widetilde{X})_q$ with the weakest one
among the topologies in which all the linear functionals $f \mapsto
\psi_{i_1i_2j_1j_2}(a)$ with $a \in q^{-2{\Bbb Z}_+}$, are continuous. The
substitution (\ref{subst}) converts $q^{-2{\Bbb Z}_+}$ into $q^{2N}\cdot
q^{-2{\Bbb Z}_+}$. After a formal passage to a limit as $N \to \infty$ we get
$q^{2{\Bbb Z}}$ instead of $q^{-2{\Bbb Z}_+}$. Equip the vector space ${\rm
Pol}(\widetilde{\Xi})_q$ with the weakest among the topologies in which all
the linear functionals $f \mapsto \psi_{i_1i_2j_1j_2}(a)$ with $a \in
q^{2{\Bbb Z}}$, $(i_1,i_2,j_1,j_2)\in M^{\prime \prime}$, are continuous. A
completion $D(\widetilde{\Xi})_q^\prime$ of this topological vector space is
the space of formal series (\ref{exp2}) whose coefficients are functions
$\psi_{i_1i_2j_1j_2}(x)$ on $q^{2{\Bbb Z}}$. The elements $f \in
D(\widetilde{\Xi})_q^\prime$ will be called distributions on the quantum
cone $\widetilde{\Xi}$. The distributions with
$\#\{(i_1,i_2,j_1,j_2,m)|\;\psi_{i_1i_2j_1j_2}(q^{-2m})\ne 0 \}<\infty$
constitute the space $D(\widetilde{\Xi})_q$ of finite functions.

 Denote by $D(\Xi)_q$, $D(\Xi)_q^\prime$ the subspaces of invariants of the
one-parameter operator groups in $D(\widetilde{\Xi})_q$,
$D(\widetilde{\Xi})_q^\prime$ which are extensions by a continuity of $i_
\varphi$, $\varphi \in{\Bbb R}/(2 \pi{\Bbb Z})$.

\begin{proposition} The structure of a covariant $*$-algebra can be
transferred by a continuity from ${\rm Pol}(\widetilde{\Xi})_q$ onto
$D(\widetilde{\Xi})_q$. The structure of a covariant
$D(\widetilde{\Xi})_q$-bimodule can be transferred by a continuity from
$D(\widetilde{\Xi})_q$ onto $D(\widetilde{\Xi})_q^\prime$.
\end{proposition}

\smallskip

 {\bf Proof.} For any polynomial $\psi$ of one indeterminate one has the
following equalities between the elements of the covariant $*$-algebra ${\rm
Pol}(\widetilde{\Xi})_q$:
$$t_{11}\psi(x)=\psi(q^2x)t_{11},\qquad t_{12}\psi(x)=\psi(x)t_{12},$$
$$t_{21}\psi(x)=\psi(x)t_{21},\qquad t_{22}\psi(x)=\psi(q^{-2}x)t_{22}.$$
Besides that, $t_{11}t_{22}=-x$, $t_{22}t_{11}=-q^{-2}x$. It follows that
the structure of an algebra in $D(\widetilde{\Xi})_q$ and the structure of a
$D(\widetilde{\Xi})_q$-bimodule in the space $D(\widetilde{\Xi})_q^\prime$
can be produced via extending by a continuity. It is easy to prove that for
any polynomial $\psi$
\begin{equation}\label{psi-mod}
\left \{\begin{array}{ccl}
X^+\psi(x)&=&-q^{1/2}t_{11}\cdot(\psi(x)-\psi(xq^{-2}))/(x-xq^{-2})\cdot
t_{21},\\
X^-\psi(x)&=&-q^{1/2}t_{12}\cdot(\psi(x)-\psi(q^{-2}x))/(x-xq^{-2})\cdot
t_{22},\\ H \psi(x)&=&0.\end{array}\right.
\end{equation}
(\ref{psi-mod}) implies that the structures of a covariant algebra in
$D(\widetilde{\Xi})_q$ and that of a covariant
$D(\widetilde{\Xi})_q$-bimodule in $D(\widetilde{\Xi})_q^\prime$ can be
produced via extending by a continuity. The rest of statements of
proposition 4.1 follow from (\ref{inv}) and $\psi(x)^*=\overline{\psi}(x)$.
\hfill $\Box$

\medskip

 Remind that the linear operator $i_ \varphi:{\rm Pol}(\widetilde{\Xi})_q
\to{\rm Pol}(\widetilde{\Xi})_q$, $\varphi \in{\Bbb R}/(2 \pi{\Bbb Z})$ is
an automorphism of covariant $*$-algebras. Hence the vector space $D(\Xi)_q
\subset D(\widetilde{\Xi})_q$ is a covariant $*$-algebra. One can prove in a
similar way that the vector space $D(\Xi)_q^\prime \subset
D(\widetilde{\Xi})_q^\prime$ are covariant $D(\Xi)_q$-bimodules.

Consider the linear functional $\nu:D(\widetilde{\Xi})_q \to{\Bbb C}$,
\begin{equation}\label{nu-i}
\int \limits_{\widetilde{\Xi}_q}fd \nu=(1-q^2)\sum_{n=-\infty}^\infty
\psi_{0000}(q^{-2m})q^{-2m},
\end{equation}
determined by the coefficient $\psi_{0000}$ in the decomposition
(\ref{exp3}).

\medskip

\begin{proposition} The linear functional (\ref{nu-i}) is an invariant
integral.
\end{proposition}

\smallskip

 {\bf Proof.} Let $\psi(x)$ be a function with finite carrier on $q^{2{\Bbb
Z}}$. One can easily deduce from the covariance of the algebra
$D(\widetilde{\Xi})_q$ and (\ref{pol-mod}), (\ref{psi-mod}) that
$$\int \limits_{\widetilde{\Xi}_q}H \psi(x)d \nu=0,\qquad \int
\limits_{\widetilde{\Xi}_q}X^+(t_{12}\psi(x)t_{22})d \nu=0.$$
Hence $\int \limits_{\widetilde{\Xi}_q}Hfd \nu=\int
\limits_{\widetilde{\Xi}_q}X^+fd \nu=0$ for all $f \in
D(\widetilde{\Xi})_q$. Thus the linear functional (\ref{nu-i}) is a $U_q
{\frak b}_+$-invariant integral. Now it remains to use the realness of this
functional. \hfill $\Box$

\medskip

 Let $l \in{\Bbb C}$. A distribution $f \in D(\widetilde{\Xi})_q^\prime$ is
called homogeneous of degree $l$ if all the coefficients
$\psi_{i_1i_2j_1j_2}(x)$ in its decomposition (\ref{exp3}) are
homogeneous:
$$\psi_{i_1i_2j_1j_2}(q^2x)=q^{2l-(i_1+i_2+j_1+j_2)}\psi_{i_1i_2j_1j_2}(x).$$
Equivalently, let $\alpha$ be an automorphism of ${\rm Pol}(\Xi)_q$ given by
$\alpha(t_{ij})=qt_{ij}$, $i,j=1,2$, and $\overline{\alpha}$ its extension
by a continuity onto $D(\widetilde{\Xi})_q^\prime$. An element $\psi \in
D(\widetilde{\Xi})_q^\prime$ is homogeneous of degree $l$ iff
$\overline{\alpha}(\psi)=q^{2l} \psi$.

 It follows from the definitions that the action of $H$, $X^+$, $X^-$ in
$D(\widetilde{\Xi})_q^\prime$ preserves the homogeneity degree  of a
distribution. Hence, the vector space $F^{(l)}$ of distributions of
homogeneity degree $l$ is a $U_q \frak{sl}_2$-module.

 Consider the linear functional
\begin{equation}\label{eta}
\eta:F^{(-1)}\to{\Bbb C};\qquad \int \limits_{\widetilde{\Xi}_q}fd
\eta=\psi_{0000}(1),
\end{equation}
determined by the coefficient $\psi_{0000}$ in the expansion (\ref{exp3}).

\medskip

\begin{proposition} The linear functional (\ref{eta}) is an invariant
integral.
\end{proposition}

\smallskip

 {\bf Proof}\ is completely similar to that of the previous proposition and
reduces to a verification of the relations
$$\int \limits_{\widetilde{\Xi}_q}Hx^{-1}d \eta=0,\qquad \int
\limits_{\widetilde{\Xi}_q}X^+(t_{12}x^{-2}t_{22})d \eta=0.$$
The first relation is obvious, and the second one follows from
$x=-t_{11}t_{22}$, $x=-qt_{12}t_{21}$, $X^+(t_{11})=X^+(t_{21})=0$. \hfill
$\Box$

 To conclude let us note that for any function on $q^{2{\Bbb Z}}$ one
has $\psi \in D(\widetilde{\Xi})_q^\prime$,
\begin{equation}\label{omega}
\Omega \psi(x)=Dx^2D\psi(x)
\end{equation}
 In the special case $\psi(x)=x^{l} $,  $l \in \Bbb C$ the next
equality from (\ref{omega}) follows:
\begin{equation}\label{omega2}
\Omega x^{l}=\frac{sh((l+1)h/2) \cdot sh(lh/2)}{sh^2(h/2)}x^{l}
\end{equation}
 One can prove the equality (\ref{omega}) by the substitution
of variables (\ref{subst}) and by the passage to a limit $N \to +\infty$
in the equality (5.7) of \cite{SSV2}.
\bigskip \setcounter{equation}{0}

\section{Cartesian products}

 Consider a Hopf algebra $A$ and $A$-module (covariant) Hopf unital algebras
$F_1$, $F_2$. Let $\pi$ be the representation of the algebra $F_1^{\rm
op}\otimes F_2$ in the vector space $F_1 \otimes F_2$ given by
$$\pi(f_1 \otimes f_2)\psi_1 \otimes \psi_2=\psi_1f_1 \otimes f_2
\psi_2;\qquad f_i,\psi_i \in F_i,\quad i=1,2.$$

\medskip

\begin{proposition} Let $K \in F_1^{\rm op}\otimes F_2$. The linear
operator $\pi(K)$ is a morphism of $A$-modules iff $K$ is an invariant:
\begin{equation}\label{inv-i}
\forall a \in A \qquad a \cdot K=\varepsilon(a) \cdot K.
\end{equation}
\end{proposition}

\smallskip

 {\bf Proof.} It follows from (\ref{inv-i}) that the linear operator $i:F_1
\otimes F_2 \to F_1^{\otimes 2}\otimes F_2^{\otimes 2}$; $i:\psi_1 \otimes
\psi_2 \mapsto \psi_1 \otimes K \otimes \psi_2$ is a morphism of
$A$-modules. Note that $\pi(K)=m_1 \otimes m_2 \cdot i$ with $m_j:F_j
\otimes F_j \to F_j$, $j=1,2$, being the multiplication in $F_j$. Hence
$\pi(K)$ is a morphism of $A$-modules. Conversely, if $\pi(K)$ is a morphism
of $A$-modules, then the elements $1 \otimes 1$ and $K=\pi(K)1 \otimes 1$ of
$F_1^{\rm op}\otimes F_2$ are invariants.\hfill $\Box$

\medskip

 We present below an evident corollary of proposition 5.1 which justifies
the algebra structure in $F_1 \otimes F_2$ introduced in \cite{SSV1}.

\begin{corollary} The invariants form a subalgebra of $F_1^{\rm op}\otimes
F_2$.  \end{corollary}

\medskip

 Let $\nu_2(f)=\int fd \nu_2$ be an invariant integral $F_2 \to{\Bbb C}$. An
importance of invariants $K \in F_1 \otimes F_2$ is due to the fact
that the associated integral operators $f \mapsto{\rm id}\otimes \nu_2(K(1
\otimes f))$ are morphisms of $A$-modules (see \cite{SSV2}).

 The previous sections contain the constructions of q-analogues for
$SU(1,1)$-spaces $\widetilde{X}$, $X$, $\widetilde{\Xi}$, $\Xi$. Let $Y \in
\{\widetilde{X},X,\widetilde{\Xi},\Xi \}$. It is easy to show that the
structure of a covariant $*$-algebra is extendable by a continuity from
${\rm Pol}(Y)_q$ onto ${\rm Fun}(Y)_q \stackrel{\rm def}{=}{\rm
Pol}(Y)_q+D(Y)_q$, and the structure of a covariant bimodule is extendable
from ${\rm Fun}(Y)_q$ onto $D(Y)_q^\prime$.

 Let $Y_1,Y_2 \in \{\widetilde{X},X,\widetilde{\Xi},\Xi \}$. Introduce the
notation
$${\rm Pol}(Y_1 \times Y_2)_q \stackrel{\rm def}{=}{\rm Pol}(Y_1)_q^{\rm
op}\otimes{\rm Pol}(Y_2)_q,$$
$$D(Y_1 \times Y_2)_q \stackrel{\rm def}{=}D(Y_1)_q^{\rm op}\otimes
D(Y_2)_q$$
for `algebras of functions on Cartesian products of quantum
$SU(1,1)$-spaces'. Note that for any of the above algebras the
multiplication is not a morphism of $U_q \frak{sl}_2$-modules, that is, the
algebras are not covariant. However, by corollary 5.2, the subspaces of
invariants are subalgebras.

 Associate to each pair of continuous linear functionals $l_j:{\rm
Pol}(Y_j)\to{\Bbb C}$, $j=1,2$, a linear functional $l_1 \otimes l_2:{\rm
Pol}(Y_1 \times Y_2)_q \to {\Bbb C}$. Equip the vector space ${\rm Pol}(Y_1
\times Y_2)_q$ with the weakest topology under which all those linear
functionals are continuous.

 The completion of the Hausdorff topological space ${\rm Pol}(Y_1 \times
Y_2)_q$ will be denoted by $D(Y_1 \times Y_2)_q^\prime$. In the case
$Y_1,Y_2 \in \{\widetilde{X},\widetilde{\Xi}\}$, $D(Y_1 \times
Y_2)_q^\prime$ may be identified with the space of formal series
\begin{equation}\label{exp4}
f=\sum_{(i_1^\prime,i_2^\prime,j_1^\prime,j_2^\prime,i_1^{\prime
\prime},i_2^{\prime \prime},j_1^{\prime \prime},j_2^{\prime
\prime})\in M^{\prime \prime}\times M^{\prime
\prime}}t_{11}^{*j_1^\prime}t_{12}^{*j_2^\prime}\otimes t_{11}^{i_1^{\prime
\prime}}t_{12}^{i_2^{\prime \prime}}\cdot \psi_{i_1^\prime i_2^\prime
j_1^\prime j_2^\prime i_1^{\prime \prime}i_2^{\prime \prime}j_1^{\prime
\prime}j_2^{\prime \prime}}(x^\prime,x^{\prime \prime}) \cdot
t_{11}^{i_2^\prime}t_{12}^{i_1^\prime}\otimes t_{11}^{*j_2^{\prime
\prime}}t_{12}^{*j_1^{\prime \prime}}.
\end{equation}
The coefficients of these formal series are functions on the Cartesian
product of progressions ($q^{-2{\Bbb Z}_+}$ or $q^{2{\Bbb Z}}$). the
topology in $D(Y_1 \times Y_2)_q^\prime$ is that of pointwise convergence of
the coefficients. It is easy to prove that the structure of $U_q
\frak{sl}_2$-module and that of ${\rm Fun}(Y_1 \times Y_2)_q$-bimodule are
extendable by a continuity from the dense linear subspace ${\rm Fun}(Y_1
\times Y_2)_q \subset D(Y_1 \times Y_2)_q^\prime$ onto the entire space
$D(Y_1 \times Y_2)_q^\prime$.

 The following statement justifies calling the bimodules $D(Y_j)_q^\prime$,
$j=1,2$, and $D(Y_1 \times Y_2)_q^\prime$ the bimodules of distributions on
the quantum $SU(1,1)$-spaces $Y_j$, $j=1,2$, and $Y_1 \times Y_2$.

 Let $\nu_j:D(Y_j)_q \to{\Bbb C}$, $\nu_j:f \mapsto \displaystyle \int
\limits_{Y_{jq}}fd \nu_j$, $j=1,2$, be invariant integrals. We follow the
conventions of \cite{SSV1} in putting in the integrands into braces the
products of elements of the algebra $F_1^{\rm op}\otimes F_2$.

\medskip

\begin{proposition} The bilinear form
\begin{equation}
D(Y_1 \times Y_2)_q \times D(Y_1 \times Y_2)_q \to{\Bbb C};\qquad \int
\limits_{Y_{1q}}\int \limits_{Y_{2q}}\{\psi_1 \cdot \psi_2 \}d \nu_{1}d
\nu_{2}
\end{equation}
is extendable by a continuity up to a bilinear form $D(Y_1 \times
Y_2)_q^\prime \times D(Y_1 \times Y_2)_q \to{\Bbb C}$. The associated map
from $D(Y_1 \times Y_2)_q^\prime$ into the dual to $D(Y_1 \times Y_2)_q$ is
an isomorphism.
\end{proposition}

\smallskip

 {\bf Proof.} It suffices to apply the description of the spaces $D(Y_1
\times Y_2)_q^\prime$ and $D(Y_1 \times Y_2)_q$ and the invariant integral
$f \mapsto \displaystyle \int \limits_{Y_{1q}}\int \limits_{Y_{2q}}fd \nu_{1}
 d \nu_{2}$
in terms of the expansion (\ref{exp4}). \hfill $\Box$

\medskip

 Proposition 5.3 implies

\medskip

\begin{corollary} There exists a unique antilinear map $\#:D(Y_1 \times
Y_2)_q^\prime \to D(Y_1 \times Y_2)_q^\prime$ such that for all $f \in
D(Y_2)_q$,  $K \in D(Y_1 \times Y_2)_q^\prime$ one has
$$\int \limits_{Y_{2q}}\{K^\# \cdot 1 \otimes f
\}d \nu_{2} =\overline{\int \limits_{Y_{2q}}\{K \cdot 1 \otimes f^*\}d
 \nu_{2}}$$.
\end{corollary}

\medskip

 Obviously, $K^{\#\#}=K$. Thus the real kernels $K=K^\#$ generate
  the `real' integral operators.

\bigskip \setcounter{equation}{0}

\section{Examples of invariant kernels}

 Let $Y_1,Y_2 \in \{\widetilde{X},\widetilde{\Xi}\}$. Consider the algebra
${\rm Pol}(Y_1 \times Y_2)_q$. We omit in the sequel the tensor product sign
$\otimes$ while working with kernels of integral operators. To avoid
misunderstanding, we introduce the notation $t_{ij},\tau_{ij}$ for the
generators $t_{ij}\otimes 1,1 \otimes t_{ij}$, $i,j=1,2$.

 The definitions imply

\medskip

\begin{proposition} The following elements of ${\rm Pol}(Y_1 \times Y_2)_q$,
$Y_1,Y_2 \in \{\widetilde{X},\widetilde{\Xi}\}$, are invariants:
$$k_{11}=t_{11}\tau_{22}-qt_{12}\tau_{21};\qquad
k_{12}=-q^{-1}t_{11}\tau_{12}+t_{12}\tau_{11};$$
$$k_{21}=t_{21}\tau_{22}-qt_{22}\tau_{21};\qquad
k_{22}=-q^{-1}t_{21}\tau_{12}+t_{22}\tau_{11}.$$
\end{proposition}

\medskip

Now prove

\medskip

\begin{proposition} The following commutation relations are valid in ${\rm
Pol}(Y_1 \times Y_2)_q$, $Y_1,Y_2 \in \{\widetilde{X},\widetilde{\Xi}\}$:
\begin{equation}
\left \{\begin{array}{lcl}k_{11}k_{12}=q^{-1}k_{12}k_{11}&,&
k_{21}k_{22}=q^{-1}k_{22}k_{21},\\ k_{11}k_{21}=q^{-1}k_{21}k_{11}&,&
k_{12}k_{22}=q^{-1}k_{22}k_{12},\\ k_{12}k_{21}=k_{21}k_{12}&,&
k_{11}k_{22}-k_{22}k_{11}=(q^{-1}-q)k_{12}k_{21}\\
k_{22}k_{11}-qk_{12}k_{21}=1 \qquad( Y_2=\widetilde{X})\\
k_{22}k_{11}-qk_{12}k_{21}=0 \qquad( Y_2=\widetilde{\Xi})
&.&\end{array}\right.
\end{equation}
\end{proposition}

\smallskip

 {\bf Proof.} We start with the special case $Y_1=Y_2=\widetilde{X}$. One
has ${\rm Pol}(Y_1 \times Y_2)_q \simeq{\Bbb C}[SL_2]_q^{\rm op}\otimes{\Bbb
C}[SL_2]_q$. Consider the map $\pi:{\Bbb C}[SL_2]_q \to{\Bbb C}[SL_2]_q^{\rm
op}\otimes{\Bbb C}[SL_2]_q$; $\pi: f \mapsto ({\rm id}\otimes S)\Delta(f)$,
with $\Delta$ being the comultiplication and $S$ the antipode of the Hopf
algebra ${\Bbb C}[SL_2]_q$:
$$\Delta(t_{ij})=\sum_kt_{ik}\otimes t_{kj},\quad S(t_{11})=t_{22},\quad
S(t_{22})=t_{11},\qquad S(t_{12})=-q^{-1}t_{12},\quad S(t_{21})=-qt_{21}.$$
The relations to be proved are provided by $\pi$ being an antihomomorphism,
and by the relations $\pi(t_{ij})=k_{ij}$, $i,j=1,2$.

 The general case reduces to the special case $Y_1=Y_2=\widetilde{X}$ via
introducing a change of generators and a passage to the limit $N \to \infty$
as in section 4. \hfill $\Box$

\medskip

 Remind that the vector space $D(Y_1 \times Y_2)_q^\prime$, $Y_1,Y_2 \in
\{\widetilde{X},\widetilde{\Xi}\}$, is equipped with the involution $\#$.

\medskip

\begin{proposition}
$$t_{11}^\#=-t_{22},\qquad t_{21}^\#=-q^{-1}t_{12},$$
\begin{equation}\label{inv-ii}
\tau_{11}^\#=-q^{-2}\tau_{22},\qquad \tau_{21}^\#=-q^{-1}\tau_{12}.
\end{equation}
\end{proposition}

\smallskip

 {\bf Proof.} Just as in the proof of proposition 6.2, we restrict ourselves
to the special case $Y_1=Y_2=\widetilde{X}$. The first two relations follow
from (\ref{inv}).

 Now to prove (\ref{inv-ii}), consider the linear functional
$l:D(\widetilde{X})_q \to{\Bbb C}$, $l:f \mapsto \displaystyle \int
\limits_{\widetilde{X}_q}f \cdot x^{-1}d \nu$. It follows from the relation
(3.3) of \cite{SSV2} that
\begin{equation}\label{l}
\forall f_1,f_2 \in D(\widetilde{X})_q \qquad l(f_1f_2)=l(f_2f_1)
\end{equation}
since $\displaystyle \int \limits_{\widetilde{X}_q}f \cdot x^{-1}d \nu={\rm
Tr}\,\widetilde{T}(f)$, with $\widetilde{T}$ being the representation of the
algebra of functions described in that work. The relations (\ref{inv-ii})
follow from
\begin{equation}\label{inv-iii}
\tau_{ij}^\#=\xi^{-1}\tau_{ij}^*\xi,\qquad i,j=1,2 \qquad with \qquad
\xi=\tau_{12} \tau_{12}^*
\end{equation}
 Finally, (\ref{inv-iii}) follow
from (\ref{l}), the realness of the invariant integral
$\nu:D(\widetilde{X})_q \to{\Bbb C}$, and the definition of the involution
$\#$. \hfill $\Box$

\medskip

 Note that, by a virtue of (\ref{l}), the linear functional $l$ may be
treated as a quantum analogue of the integral with respect to the Liouville
measure associated to the standard symplectic structure on
$\widetilde{X}\subset SL_2({\Bbb C})$.

 Proposition 6.3 allows one to produce explicit formulae for an involution
in the algebra of invariant polynomial kernels. Specifically, one has

\medskip

\begin{corollary} $k_{11}^\#=q^{2}k_{22},\; k_{12}^\#=q^{-1}k_{21}$.
\end{corollary}

\medskip

 Our immediate purpose is to study the invariants $k_{22}^lk_{11}^l\in D(Y_1
\times Y_2)_q^\prime$.

 Remind the standard notation $(t;q)_n=\prod \limits_{j=0}^{n-1}(1-tq^j)$.
It is easy to prove (see \cite[relation(1.3.2)]{GR}) that
\begin{equation}\label{tqn}
(t;q)_n=\sum_{j=0}^n \frac{(q^{-n};q)_j}{(q;q)_j}q^{j(n+1)}t^j.
\end{equation}

 Just as in section 3, introduce the linear subspace $F \subset D(Y_1
\times Y_2)_q^\prime$ of finite sums of the form (\ref{exp4}). The structure
of algebra is transferred by a continuity from ${\rm Pol}(Y_1 \times
Y_2)_q$ onto $F$. Impose the notation for some special elements of $F$:
$$z=qt_{12}^{-1}t_{11},\quad z^*=t_{22}t_{21}^{-1},\quad \zeta=q
\tau_{11}\tau_{12}^{-1},\quad \zeta^*=\tau_{21}^{-1}\tau_{22},\quad
x=t_{12}^*t_{12},\quad \xi=\tau_{12}\tau_{12}^*.$$

 Use the commutation relations
$$t_{12}\tau_{21}(z \zeta^*)=q^2(z \zeta^*)t_{12}\tau_{21};\qquad (\zeta
z^*)t_{12}\tau_{21}=q^2t_{21}\tau_{12}(\zeta z^*)$$
to prove following relations in $F$.

 Let $l \in{\Bbb Z}_+$. Then
\begin{equation}\label{kl1}
{k_{11}^l=(-qt_{12}\tau_{21}(1-q^{-2}z
\zeta^*))^l=(-qt_{12}\tau_{21})^l(q^{-2}z \zeta^*;q^{-2})_l,\atop
k_{22}^l=((1-z^*\zeta)(-q^{-1}t_{21}\tau_{12}))^l=
(z^*\zeta;q^{-2})_l(-q^{-1}t_{21}\tau_{12})^l.}
\end{equation}
Hence,
$$k_{22}^lk_{11}^l=
(z^*\zeta;q^{-2})_l(t_{12}t_{21})^l(\tau_{12}\tau_{21})^l(q^{-2}z
\zeta^*;q^{-2})_l;$$
$$k_{22}^lk_{11}^l=q^{-2l}\xi^l(q^2z^*\zeta;q^2)_l(z \zeta^*;q^2)_lx^l.$$

 By a virtue of (\ref{tqn}) one has
\begin{equation}\label{l1}
{(q^2z^*\zeta;q^2)_l=\displaystyle \sum_{n=0}^\infty
\frac{(q^{-2l};q^2)_n}{(q^2;q^2)_n}q^{2(l+1)n}(z^*\zeta)^n,\atop(z
\zeta^*;q^2)_l=\displaystyle \sum_{n=0}^\infty
\frac{(q^{-2l};q^2)_n}{(q^2;q^2)_n}q^{2ln}(z \zeta^*)^n.}
\end{equation}
Thus we have proved

\medskip

\begin{lemma} For all $l \in{\Bbb Z}_+$,
\begin{equation}\label{kl2}
k_{22}^lk_{11}^l=q^{-2l}\xi^l \sum_{j=0}^\infty
\frac{(q^{-2l};q^2)_j}{(q^2;q^2)_j}(q^{2(l+1)}z^*\zeta)^j \sum_{m=0}^\infty
\frac{(q^{-2l};q^2)_m}{(q^2;q^2)_m}(q^{2l}z \zeta^*)^m x^l.
\end{equation}
\end{lemma}

\medskip

 {\sc Remark 6.6.} Show that the right hand side of (\ref{kl2}) is a
generalized kernel, that is, it could be written in the form (\ref{exp4}).
$$k_{22}^lk_{11}^l=q^{-2l}\sum_{j=0}^\infty \sum_{m=0}^\infty
\frac{(q^{-2l};q^2)_j}{(q^2;q^2)_j}
\frac{(q^{-2l};q^2)_m}{(q^2;q^2)_m}q^{2j(l+1)+2ml}z^{*j}z^mx^l \xi^l \zeta^j
\zeta^{*m}=$$
$$=q^{-2l}\sum_{0 \le j<m<\infty}\frac{(q^{-2l};q^2)_j}{(q^2;q^2)_j}
\frac{(q^{-2l};q^2)_m}{(q^2;q^2)_m}q^{2j(l+1)+2ml}(z^{*j}z^j)z^{m-j}x^l
\xi^l(\zeta^j \zeta^{*j})\zeta^{*(m-j)}+$$
$$+q^{-2l}\sum_{j=0}^\infty
\frac{(q^{-2l};q^2)_j^2}{(q^2;q^2)_j^2}q^{4jl+2j}(z^{*j}z^j)x^l
\xi^l(\zeta^j \zeta^{*j})+$$
$$+q^{-2l}\sum_{0 \le j<m<\infty}\frac{(q^{-2l};q^2)_j}{(q^2;q^2)_j}
\frac{(q^{-2l};q^2)_m}{(q^2;q^2)_m}q^{2j(l+1)+2ml}z^{*(j-m)}(z^{*m}z^m)x^l
\xi^l \zeta^{j-m}(\zeta^m \zeta^{*m}).$$

\medskip \stepcounter{theorem}

 Now we are in a position to make up an analytic continuation of the
distribution $k_{22}^lk_{11}^l$ in $l$. Consider a vector-function
$f(\lambda)$ of a complex variable $\lambda$ with values at $D(Y_1 \times
Y_2)_q^\prime$, $Y_1,Y_2 \in \{\widetilde{X},\widetilde{\Xi}\}$. This
vector-function will be called polynomial if for any finite function $\psi
\in D(Y_1 \times Y_2)_q$ the integral $I(\lambda)=\int \limits_{Y_{1q}}\int
\limits_{Y_{2q}}f(\lambda)\psi d \nu_{1}d \nu_{2}$ is a polynomial of
$\lambda,\lambda^{-1}$, that is, $I(\lambda)\in{\Bbb
C}[\lambda,\lambda^{-1}]$.

\medskip

\begin{proposition} Consider any one among the three spaces of distributions
$D(\widetilde{X}\times \widetilde{X})_q^\prime$, $D(\widetilde{X}\times
\widetilde{\Xi})_q^\prime$, $D(\widetilde{\Xi}\times
\widetilde{X})_q^\prime$. There exists a unique polynomial vector-function
$f(\lambda)$ with values at a selected space of distributions such that
$f(q^{2l})=k_{22}^lk_{11}^l$ for all $l \in{\Bbb Z}_+$.
\end{proposition}

\smallskip

 {\bf Proof.} The uniqueness of a rational function $f(\lambda)$ with given
(fixed) values in $\lambda=q^{2l}$, $l \in{\Bbb Z}_+$ is evident. Consider a
finite function $\psi$ and replace $k_{22}^lk_{11}^l$ in the integral $\int
\limits_{Y_{1q}}\int \limits_{Y_{2q}}f(\lambda)\psi d \nu_{1}d \nu_{2}$ by the
right hand side of (\ref{kl2}). If either $Y_1=\widetilde{X}$ or
$Y_2=\widetilde{X}$, then only finitely many terms of (\ref{kl2}) contribute
to the above integral. So, what remains is to note that each term is a
polynomial of $q^{2l}$, $q^{-2l}$. \hfill $\Box$

\medskip

 {\sc Remark 6.8.} A distribution $f(q^{2l})$ whose existence and uniqueness
was proved in proposition 6.6, will be also denoted by $k_{22}^lk_{11}^l$.
By corollary 5.2, it is invariant for all $l \in{\Bbb Z}_+$, and hence for
all $l \in{\Bbb C}$.

\medskip \stepcounter{theorem}

To conclude, we prove the following

\begin{proposition} For any of the spaces $D(\widetilde{X}\times
\widetilde{X})_q^\prime$, $D(\widetilde{X}\times \widetilde{\Xi})_q^\prime$,
$D(\widetilde{\Xi}\times \widetilde{X})_q^\prime$, there exists a unique
polynomial vector-function $f(\lambda)$ such that for all $l \in{\Bbb Z}_+$
$$k_{11}^{-l}=(-qt_{12}\tau_{21})^{-l}f(q^{2l}).$$
\end{proposition}

\smallskip

 {\bf Proof.} The uniqueness is obvious. The existence follows from
(\ref{kl1}), (\ref{l1}), and the q-binomial theorem (see \cite{GR}).
\begin{equation}
f=(z \zeta^*;q^2)_l^{-1}=\sum_{n=0}^\infty
\frac{(q^{2l};q^2)_n}{(q^2;q^2)_n}(z \zeta^*)^n.
\end{equation}

\bigskip \setcounter{equation}{0}

\begin{appendix}
\section*{Appendix. On involution in ${\rm Pol}(\widetilde{X})_q$}

 The standard procedure of constructing an involution in ${\Bbb C}[SL_2]_q$
is to use the involution $*:U_q \frak{su}(1,1)\to U_q \frak{su}(1,1)$, the
antilinear map
$$\#:\xi \mapsto(S(\xi))^*,\qquad \xi \in U_q \frak{su}(1,1),$$
and the duality argument
$$\forall \xi \in U_q \frak{su}(1,1),\;f \in{\Bbb C}[SL_2]_q \qquad
f^*(\xi)\stackrel{\rm def}{=}\overline{f(\xi^*)}.$$
One gets in this way a Hopf $*$-algebra ${\Bbb C}[SU(1,1)]_q$ of
regular functions on the quantum group $SU(1,1)$.

 In order to obtain the involution involved in the definition of the
principal quantum homogeneous space, one can modify the above procedure. Let
$w_q \in{\Bbb C}[SL_2]_q^*$ be the element of the quantum Weyl group
\cite{CP}. Consider the linear subspace $L \stackrel{\rm def}{=}w_q^{-1}U_q
\frak{sl}_2=U_q \frak{sl}_2w_q^{-1}$ (The last equality follows from
$$w_q \cdot K^{\pm 1}\cdot w_q^{-1}=K^{\pm 1},\qquad w_q \cdot E \cdot
w_q^{-1}=-q^{-1}F,\qquad w_q \cdot F \cdot w_q^{-1}=-qE \eqno(A.1)$$
(see \cite{KR})). The `involution' $\#:U_q \frak{su}(1,1)\to U_q
\frak{su}(1,1)$ is to be replaced by the `involution'
$$\#:L \to L,\qquad \#:w_q^{-1}\xi \mapsto w_q^{-1}(S(\xi))^*.$$
Finally, let $*$ be such an antilinear operator in ${\Bbb C}[SL_2]_q$ that
$$\forall \xi \in L,\;f \in{\Bbb C}[SL_2]_q \qquad
f^*(\xi)=\overline{f(\xi^\#)}.$$

 Evidently, $f^{**}=f$ for all $f \in{\Bbb C}[SL_2]_q$.

 Prove that $*$ is an antihomomorphism of ${\Bbb C}[SL_2]_q$. This follows
from
$$\Delta w_q^{-1}=w_q^{-1}\otimes w_q^{-1}\cdot R,\qquad
R^{*\otimes*}=R_{21},\qquad S \otimes S(R)=R,$$
with $R$ being a universal R-matrix, and $R_{21}$ is derived from $R$ by a
permutation of tensor multiples (cf. \cite{SV}). It is worthwhile to note
that for all $f \in{\Bbb C}[SL_2]_q$ and all $\xi \in U_q \frak{su}(1,1)$
one has $(\xi f)^*=(S(\xi))^*f^*$. Thus we get a covariant $*$-algebra.
Obviously, the linear functional $w_q$ is real:
$f^*(w_q)=\overline{f(w_q)}$ for all $f \in{\Bbb C}[SL_2]_q$.

 To conclude, let us prove that the above involution coincides with that in
${\rm Pol}(\widetilde{X})_q$, i.e. $t_{11}^*=-t_{22}$, $t_{12}^*=-qt_{21}$.
It follows from the definition of the involution that the linear span of
$t_{ij}\in{\Bbb C}[SL_2]_q$, $i,j=1,2$, is an invariant subspace for $*$.
What remains is to apply the realness of the functional $w_q \in{\Bbb
C}[SL_2]_q^*$ and the relations
$$\pmatrix{t_{11}(w_q) & t_{12}(w_q)\cr t_{21}(w_q) & t_{22}(w_q)}={\rm
const}(q)\pmatrix{0 & -q \cr 1 & 0},$$ with ${\rm const}(q)\to 1$ as $q \to
1$ (see \cite{KR}).

\end{appendix}


\bigskip

\begin{thebibliography}{99}

\bibitem{CP} V. Chari, A. Pressley. A Guide to Quantum Groups, Cambridge
Univ. Press, 1995.

\bibitem{GR} G. Gasper, M. Rahman. Basic Hypergeometric Series, Cambridge
University Press, Cambridge, 1990.

\bibitem{KR} A. N. Kirillov, Ya. S. Reshetikhin, {\it q-Weyl group and a
multiplicative formula for universal R-matrices}, Commun. Math. Phys., {\bf
134}, (1990), 421 -- 431.

\bibitem{SSV1} D. Shklyarov, S. Sinel'shchikov, L. Vaksman. On function
theory in quantum disc: integral representations, E-print: math.QA/9808015.

\bibitem{SSV2} D. Shklyarov, S. Sinel'shchikov, L. Vaksman. On function
theory in quantum disc: covariance, E-print: math.QA/9808037.

\bibitem{SV} S. Sinel'shchikov, L. Vaksman. {\it On q-analogues of Bounded
Symmetric Domains and Dolbeault Complexes},  Mathematical Physics, Analysis
and Geometry; Kluwer Academic Publishers, V.1, No.1, 1998, 75--100,
E-print: q-alg/9703005.


\bibitem{V} L. Vaksman, {\it q-analogues of the Clebsch-Gordan coefficients
and algebra of functions on the quantum group $SU(2)$}, Dokl. Acad. Sci.
USSR, {\bf 306} (1989), No 2, 269 -- 271.

\end{thebibliography}
\end{document}